\documentclass[11pt]{amsart}
\usepackage{graphicx}
\usepackage{amssymb}
\usepackage{epstopdf}
\usepackage{amsmath}
\usepackage{amscd}
\usepackage{amsfonts,amsxtra,amssymb,amscd}
\usepackage[all]{xy}
\usepackage{enumerate}
\usepackage{hyperref}

\newtheorem{theorem}{Theorem}[section]

\newtheorem{lemma}[theorem]{Lemma}

\def\lato{\vbox {\hrule width 1.5cm height2pt depth-1.5pt }}
\def\llato{\vbox {\hrule width 1.5cm height.8pt depth-.3pt }}
\def\lllato{\vbox {\vskip-1.5pt\hrule width 1.5cm height.4pt}}
\def\llllato{\vbox {\hrule width 1.5cm height.4pt}}
\def\bblato{\overset \lllato    \lllato}
\def\tlato{\overset \llllato\bblato \hskip-0.85cm >
\ \ \ \,}
\def\plato{\vbox {\hrule width 1cm height2pt depth-1.5pt }}
\def\blato{\overset\llato \llato \hskip-0.85cm <
\ \ \ \,}
\def\pblato{\overset\llato \llato \hskip-0.85cm >
\ \ \ \,}

\def\F4{  {\underset  1\circ } \lato {\underset  2\circ } \blato {\underset 
3\circ } \lato  {\underset  4\circ }}
\def\G2{  {\underset  1\circ } \tlato {\underset  2\circ }
}
\def\Dn{  {\underset 1\circ }\lato {\underset  2\circ }
\cdots \cdots \plato{\underset  {n-2}\circ }\lato {\underset  {n-1}\circ }\raise5pt  \hbox { 
\hskip-2.42cm
   \vrule width 0.5pt height1.5cm }\raise1.7cm \hbox {\hskip-2.7pt$
\overset n \circ$}   }
\def\esei{  {\underset 1\circ }\lato {\underset  2\circ } \lato {\underset 
3\circ }\lato {\underset  4\circ
}\lato {\underset  5\circ
}\raise5pt  \hbox {  \hskip-3.58cm
   \vrule width 0.5pt height1.5cm }\raise1.7cm
\hbox {\hskip-2.7pt$\overset 6 \circ$}}

\def\esette{  {\underset  1\circ }\lato {\underset  2\circ } \lato {\underset 
3\circ }\lato {\underset  4\circ
}\lato {\underset  5\circ
}\lato {\underset  6\circ
}\raise5pt  \hbox {  \hskip-3.58cm
   \vrule width 0.5pt height1.5cm }\raise1.7cm \hbox {\hskip-2.7pt$\overset 7 \circ$}}

\def\E8{  {\underset  1\circ }\lato {\underset  2\circ } \lato {\underset 
3\circ }\lato {\underset  4\circ
}\lato {\underset  5\circ
}\lato {\underset  6\circ
}\lato {\underset  7\circ
}\raise5pt  \hbox {  \hskip-3.58cm
   \vrule width 0.5pt height1.5cm }\raise1.7cm \hbox {\hskip-2.7pt$\overset 8 \circ $}}

\def\AfBn{ \hskip-5.5cm {\underset  1\circ } \lato {\underset  2\circ } \lato {\underset 
3\circ }  \cdots\cdots \cdots \pblato  {\underset  n\circ } 
\raise5pt  \hbox { 
\hskip-5.33cm
   \vrule width 0.5pt height1.5cm }\raise1.7cm \hbox {\hskip-2.7pt$
\overset 0 \circ$}  }
\def\AfCn{ {\underset  0\circ }\pblato {\underset  1\circ } \lato {\underset  2\circ } \lato {\underset 
3\circ }  \cdots\cdots \cdots \blato {\underset  n\circ } }
\def\AfF4{ {\underset  0\circ } \lato {\underset  1\circ } \lato {\underset  2\circ } \pblato {\underset 
3\circ } \lato  {\underset  4\circ }}
\def\AfG2{ {\underset  0\circ }\lato  {\underset  1\circ } \tlato {\underset  2\circ }
}
\def\AfDn{\hskip-5.5cm  {\underset 1\circ }\lato {\underset  2\circ }
\cdots \cdots \plato{\underset  {n-2}\circ }\lato {\underset  {n-1}\circ }\raise5pt  \hbox { 
\hskip-2.48cm
   \vrule width 0.5pt height1.5cm }\raise1.7cm \hbox {\hskip-2.7pt$
\overset n \circ$}\raise5pt  \hbox { 
\hskip-2.73cm
   \vrule width 0.5pt height1.5cm }\raise1.7cm \hbox {\hskip-2.7pt$
\overset 0 \circ$}  }

\def\Afesei{\hskip-5.5cm  {\underset 1\circ }\lato {\underset  2\circ } \lato {\underset 
3\circ }\lato {\underset  4\circ
}\lato {\underset  5\circ
}\raise5pt  \hbox {  \hskip-3.62cm
   \vrule width 0.5pt height1.5cm }\raise1.7cm
\hbox {\hskip-2.7pt$\overset {\ 6} \circ$}\raise1.9cm  \hbox { \hskip-9.0pt
   \vrule width 0.5pt height1.5cm }\raise3.4cm
\hbox {\hskip-2.7pt$\overset 0 \circ$}}

\def\Afesette{\hskip-4.5cm   {\underset  0\circ }\lato{\underset  1\circ }\lato {\underset  2\circ } \lato {\underset 
3\circ }\lato {\underset  4\circ
}\lato {\underset  5\circ
}\lato {\underset  6\circ
}\raise5pt  \hbox {  \hskip-5.32cm
   \vrule width 0.5pt height1.5cm }\raise1.7cm \hbox {\hskip-2.7pt$\overset 7 \circ$}}

\def\AfE8{\hskip-3.2cm{\underset  0\circ }\lato  {\underset  1\circ }\lato {\underset  2\circ } \lato {\underset 
3\circ }\lato {\underset  4\circ
}\lato {\underset  5\circ
}\lato {\underset  6\circ
}\lato {\underset  7\circ
}\raise5pt  \hbox {  \hskip-3.63cm
   \vrule width 0.5pt height1.5cm }\raise1.7cm \hbox {\hskip-2.7pt$\overset 8 \circ $}}

\title[A curious identity]{A curious identity and the volume of the root spherical simplex.
}
\author{C. De Concini}\address{Dip. Mat. Castelnuovo, Univ. di Roma La
Sapienza, Rome, Italy}\email{deconcin@mat.uniroma1.it}
\author{C. Procesi}\address{Dip. Mat. Castelnuovo, Univ. di Roma La
Sapienza, Rome, Italy}\email{procesi@mat.uniroma1.it} 
\address{Dept.~of Mathematics, University of Michigan, Ann Arbor MI 48109 USA}
\email{jrs@umich.edu}
\thanks{The authors are partially supported by the Cofin 40
\%, MIUR}

\begin{document}
\dedicatory{{ \rm with an appendix by}\vspace{10pt}\\ {\large {\textsc{john r. stembridge}}}\vspace{10pt} \\ A Guido Zappa per i suoi 90 anni
\\} 
\begin{abstract}{We show a curious identity on root systems which gives the evaluation of the volume of the spherical simpleces cut by the cone generated by simple roots. In the appendix John Stembridge gives a conceptual proof of our identity}\end{abstract}
\maketitle

  \section{Introduction}

In this note we shall consider a finite root system $R$ spanning an euclidean space $E$ of dimension $\ell$ (for all the  facts about root systems which we are going to use in this note we refer to \cite{bou}). $\ell$ is called the rank of $R$. We shall choose once and for all a set of positive roots $R^+$ and in   $R^+$ the set of simple roots $\Delta=\{\alpha_1,\ldots ,\alpha_{\ell}\}$. We shall also denote by $W$ the Weyl group of $R$ i.e. the finite group generated by the reflections with respect to the hyperplanes orthogonal to the roots in $R$. Given such a root $\alpha\in R$, we shall denote by $s_{\alpha}\in W$ the reflection with respect to the hyperplane orthogonal to $\alpha$.  Set $s_i=s_{\alpha_i}$ for each $i=1,\ldots ,\ell$ and call $S=\{s_1,\ldots ,s_{\ell}\}$ the set of simple reflections. One know that $S$ generates $W$ and that   the pair $(W,S)$ is a Coxeter group.

We know that the ring of regular  functions on $E$, invariant  under the action of $W$, is a polynomial ring generated by homogenous elements of degrees $d_1\leq d_2\leq\cdots \leq d_{\ell}$. The $d_i$'s are called the degrees. We shall also consider the sequence of exponents, $d_1-1,d_2-1,\ldots ,d_{\ell}-1$. Recall that $\prod_id_i=|W|$.

In $E$ we have the {\it affine arrangement}  of the hyperplanes orthogonal to the roots and  their translates under  the weight lattice $\Lambda$,  a locally finite configuration  invariant under the affine Weyl group  $\hat W$.  $\hat W$ is the semidirect product of $W$ and of the lattice $Q$ spanned by the roots, thought of as translation operators. 

 $\hat W$ is itself a Coxeter group. In the case in which $E$ is irreducible, its Coxeter generators are given by the reflections $\{s_0, s_1\ldots , s_{\ell})$, where the for $i\geq 1$ the $s_i$'s are the simple generators of $W$ and 
$$s_0(v)=s_{\theta}(v)+\theta$$
$\theta$ being the longest root. One knows that, for each $0\leq i\leq \ell$, the subgroup $W_i$ of $\hat W$ generated by the reflections $(s_0,\ldots ,\check s_i,\ldots ,s_{\ell})$ is finite, and it is the Weyl group of a root system $R^{(i)}$ which will be discussed presently.
 Hence we can consider the degrees $d^{(i)}_1\leq d^{(i)}_2\leq\cdots \leq d^{(i)}_{\ell}$. Our main result is
the identity (Theorem \ref{ident})
$$\sum_{i=0}^\ell \frac  {(d^{(i)}_1-1)(d^{(i)}_2-1)\cdots ( d^{(i)}_{\ell}-1)}{d^{(i)}_1d^{(i)}_2\cdots d^{(i)}_{\ell}}=1.$$
The proof is a case by case computation using the classification of irreducible root systems. It is quite desirable to give a more conceptual deduction of our identity.

In the last section we show, following a suggestion of Vinberg,  that our identity implies the following geometric identity.  Take the unit sphere $S(E)$ in $E$ and consider the spherical simplex $S(E)=C(\Delta)\cap S(E)$, $C(\Delta)$ being the cone of positive linear combinations of the simple roots. Then
$$\frac {{\rm Vol\ }S(\Delta)}{{\rm Vol\ }S(E)}= \frac  {(d_1-1)(d_2-1)\cdots ( d_{\ell}-1)}{d_1d_2\cdots d_{\ell}}.$$

We have discovered this identity while trying to understand the following fact. 

Consider the complex
space $V=E\otimes_{\mathbb R}\mathbb C$, and take the algebraic torus $T=V/Q$.  For any root $\alpha\ \in R$ the linear form 
$\alpha\check{}\ $ 
defined by $$ \alpha\check{}\ (v)=2\frac{(\alpha,v)}{(\alpha,\alpha)}$$ takes integer values on $Q$, hence we get the character $e^{2\pi\sqrt{-1} \alpha\check{}\  }$ of $T$.  

Denote its kernel  by $D_\alpha$. In our work on toric arrangements (see \cite{dp} ,\cite{He}, \cite{Le1}, \cite{Le2}, \cite{lo}) we have shown  that the Euler characteristic of the open set $\mathcal A:=T-\cup_{\alpha\in R^+}D_{\alpha}$ equals $(-1)^{\ell}|W|$. The only proof we know of this fact is via a combinatorial topological construction of Salvetti \cite{S3} \cite{dS3}. The above identity has been the result of an attempt to 
give a direct computation of this Euler characteristic.

 \subsection{The main  identity}  We are interested in the numbers
$$\nu(R)=\prod_{i=1}^{\ell}{d_i-1\over d_i}.$$

The following table gives $\nu(R)$ in the case of irreducible root systems
$$ A_n, n\geq 1,\quad  \nu_{A_n}={1\over n+1}$$
$$ B_n {\rm \ and\ } C_n, n\geq 2,\quad    \nu_{B_n}= {1\over 4^n}{\binom {2n} n}$$
$$D_n, n\geq 4, \quad  \nu_{D_n}= {n-1\over 4^{n-1}n}{\binom {2(n-1)} {n-1}}$$
$$\nu_{G_2}= {5\over 12}$$
$$\nu_{F_4}={385\over 1152}$$ 
$$\nu_{E_6}= {77\over 324}$$
$$\nu_{E_7}= {2431\over 9216}$$
$$\nu_{E_8}={30808063\over 99532800}$$

Notice that if $R$ is reducible, i.e. $R=R_1\cup R_2$ with $R_1\perp R_2$, then clearly
$$\nu(R)=\nu(R_1)\nu(R_2).$$
We normalize the scalar product so that the short roots have length $\sqrt 2$ and we denote by $D$ the Dynkin diagram of $R$.

From now on we assume that $R$ is irreducible and we denote by $\hat D$ the extended Dynkin diagram.  Set $\alpha_0=-\theta$, with  $\theta$ the highest root. 
We have a bijection between the set $\hat\Delta =\{\alpha_0,\ldots ,\alpha_\ell\}$ and the nodes of  $\hat D$. For every  $i=0\ldots ,\ell$  the diagram $D_i$ obtained from $\hat D$ removing the node corresponding to the root $\alpha_i$ (and all the edges having that node as one of the vertices) is of finite type. So  we can consider
the corresponding root system $R^{(i)}$  consisting of all roots in $R$ which are integral linear combinations of the roots $\alpha_0,\ldots ,\check\alpha_i ,\ldots ,\alpha_{\ell}$ and   the corresponding number $\nu(R^{(i)})$.  

During the proof of our result we shall need the following well known
\begin{lemma}\label{illem} The following identities hold,
 \begin{equation}4^n=\sum_{h=0}^n\binom{2h}{ h}\binom{2(n-h)}{n-h}, \ n\geq 0\end{equation}
 Furthermore, when $n\geq 2$ we have:
\begin{equation}4^{n-1}={1\over 2}\binom{2n}{ n}+\sum_{h=2}^n{h-1\over h}\binom{2(h-1)}{ h-1}\binom{2(n-h)}{ n-h} \end{equation}

  \begin{equation}4^{n-2}={n-1\over 4n}\binom {2(n-1)} {n-1}+\sum_{h=2}^{n-2}{(h-1)(n-h-1)\over h(n-h)}\binom{2(h-1)}{ h-1}\binom{2(n-1-h)}{ n-1-h}.\end{equation}
 \end{lemma}
\begin{proof} The first identity follows immediately from the following power series expansion
\begin{equation}\label{prima} {1\over \sqrt{1-4t}}=\sum_{n\geq 0}\binom{2n}{n}t^n.
\end{equation}
To see this notice that setting $$f(t)={1\over \sqrt{1-4t}}$$ we have
$${df(t)\over dt}=2f(t)^3$$
from which we deduce that
$$(1-4t){df(t)\over dt}=2f(t).$$
Writing $f(t)=\sum_{n\geq 0} a_nt^n$ we deduce that
$$\sum_{n\geq 0} na_nt^{n-1}-4\sum_{n\geq 0} na_nt^n=2\sum_{n\geq 0} a_nt^n.$$
Equating coefficients, we get $a_0=1$ and for $h\geq 0$
$$a_{h+1}={2(2h+1)\over h+1}a_h.$$
On the other hand if we set $b_h:=\binom{2h}{h}$, we get $b_0=1$ and
$$b_{h+1}={2(h+1)!\over (h+1)!(h+1)!}={2(h+1)(2h+1)\over (h+1)^2}b_h={2(2h+1)\over h+1}b_h$$
so $a_n=b_n$ and everything follows.

To see the second identity, notice that using (\ref{prima}) and integrating, we get
\begin{equation}{1\over 2}-\sum_{h\geq 1}{1\over h}\binom{2(h-1)}{ h-1}t^h={1\over 2}\sqrt{1-4t}\end{equation}
Again using (\ref{prima}) we deduce
\begin{equation}{1\over 2}+\sum_{h\geq 2}{h-1\over h}\binom{2(h-1)}{ h-1}t^h=\sum_{h\geq 1}\binom{2(h-1)}{ h-1}t^h+{1\over 2}-\sum{1\over h}\binom{2(h-1)}{ h-1}t^h=$$$$= {1-2t\over 2\sqrt{1-4t}}\end{equation}
This together with (\ref{prima}) implies that 
$${1\over 2}\binom{2n}{ n}+\sum_{h=2}^n{h-1\over h}\binom{2(h-1)}{ h-1}\binom{2(n-h)}{ n-h}$$
is the coefficient of $t^n$ in the power series expantion of
$$ {1-2t\over 2\sqrt{1-4t}} {1\over \sqrt{1-4t}}={1\over 2}+{t\over 1-4t}$$
since $n\geq 2$ the claim follows.
To see the last identity, let us remark that its left handside is the coefficient of $t^n$ in the power series expansion of the function
$$({1-2t\over 2\sqrt{1-4t}})^2={1+4t^2-4t\over 1-4t}.$$
From this everything follows.

\end{proof}
\begin{theorem}\label{ident} $\sum_{i=0}^\ell\nu( R^{(i)})=1$
\end{theorem}
\begin{proof} The proof is by a case by case computation.

Let us deal first with the exceptional cases. In order to make the computation transparent it is more convenient to multiply our sum by $|W|$, the order of the Weyl group

Case $G_2$.   In this case $|W|=12$. By looking at the extended Dynkin diagram
$$\AfG2$$
we get that
$$|W|\sum_{i=0}^2\nu ( R^{(i)})=12(\nu(G_2)+\nu(A_2)+\nu(A_1\times A_1))=5+3+4=12$$
Case $F_4$. The order of the Weyl group is $1152$. By looking at the extended Dynkin diagram
$$\AfF4$$
we get that
$$|W|\sum_{i=0}^4\nu ( R^{(i)})=|W|(\nu(F4)+\nu(A_1\times C_3)+\nu(A_2\times A_2)+\nu(A_3\times A_1)+\nu(B_4))=$$$$=385+180+128+144+315=1152.$$
Case $E_6$.   In this case $|W|=51840$. By looking at the extended Dynkin diagram
$$\Afesei$$
we get that
$$|W|\sum_{i=0}^6\nu ( R^{(i)})=|W|(3\nu(E_6)+\nu(A_2\times A_2\times A_2)+3\nu(A_1\times A_5))=$$$$=36960+1920+12960=51840$$
Case $E_7$.   In this case $|W|=2903040$. By looking at the extended Dynkin diagram
$$\Afesette$$
we get that
$$|W|\sum_{i=0}^7\nu ( R^{(i)})=|W|(2\nu(E_7)+2\nu(A_1\times D_6)+2\nu(A_2\times A_5)+\nu(A_3\times A_3\times A_1)+\nu(A_7))=$$$$=1531530+595350+322560+90720+362880=2903040$$
Case $E_8$.   In this case $|W|=696729600$. By looking at the extended Dynkin diagram
$$\AfE8$$
we get that
$$|W|\sum_{i=0}^8\nu ( R^{(i)})=|W|(\nu(E_8)+\nu(A_1\times E_7)+\nu(A_2\times E_6)+\nu(A_3\times D_5)+\nu(A_4\times A_4)+$$$$+\nu(A_5\times A_2\times A_1)+\nu(A_1\times A_7)+\nu(D_8)+\nu(A_8))=$$$$=215656441+91891800+55193600+38102400+27869184+$$$$+19353600+43545600+127702575+77414400=696729600$$
Case $A_n$. In this case each $ R^{(i)}$ is of type $A_n$. It follows that 
$$\sum_{i=0}^n\nu( R^{(i)})=(n+1)\nu_{A_n}=(n+1){1\over n+1}=1$$
Case $C_n$.   The extended Dynkin diagram is
$$\AfCn$$
we get that, denoting by $C_0$ the trivial root system and setting $C_1=A_1$,

$$\sum_{i=0}^n\nu ( R^{(i)})=\sum_{h=0}^n\nu(C_h\times C_{n-h})=\sum_{h=0}^n{1\over 4^{h+n-h}}\binom{2h}{ h}\binom{2(n-h)}{n-h}=1$$
by   Lemma \ref{illem}, part (1).

Case $B_n$.   The extended Dynkin diagram is
$$\AfBn$$
we get that, denoting by $B_0$ the trivial root system, setting $C_1=A_1$,  $D_2=A_1\times A_1$  and  $D_3=A_3$,
$$\sum_{i=0}^n\nu ( R^{(i)})=2\nu(B_n)+\sum_{h=2}^n\nu(D_h\times B_{n-h})=$$$$={2\over 4^n}\binom{2n}{ n}+\sum_{h=2}^n{h-1\over h}\binom{2(h-1)}{ h-1}\binom{2(n-h)}{ n-h}=$$$$={1\over 4^{n-1}}({1\over 2}\binom{2n}{ n}+\sum_{h=2}^n{h-1\over h}\binom{2(h-1)}{ h-1}\binom{2(n-h)}{ n-h})=1$$
by   Lemma \ref{illem}, part (2).

Case $D_n$.   The extended Dynkin diagram is
$$\AfDn$$
we get that, setting  $D_2=A_1\times A_1$  and  $D_3=A_3$,
$$\sum_{i=0}^n\nu ( R^{(i)})=4\nu(D_n)+\sum_{h=2}^{n-2}\nu(D_h\times D_{n-h})=$$$${1\over 4^{n-2}}({n-1\over 4n}\binom {2(n-1)} {n-1}+\sum_{h=2}^{n-2}{(h-1)(n-h-1)\over h(n-h)}\binom{2(h-1)}{ h-1}\binom{2(n-1-h)}{ n-1-h}) $$ which equals 1 by   Lemma \ref{illem}, part (3).
\end{proof}
\subsection{The volume of $S(\Delta)$} Recall that we have introduced the spherical simplex as the intersection of the unit sphere $S(E)$ in $E$ with the cone $C(\Delta)$ of non negative linear combinations of the simple roots $\{\alpha_1,\alpha_2\ldots ,\alpha_{\ell}\}$ for the root system $R$. Our purpose is to show
\begin{theorem}  $$\frac {{\rm Vol\ }S(\Delta)}{{\rm Vol\ }S(E)}=\nu(R)= \frac  {(d_1-1)(d_2-1)\cdots ( d_{\ell}-1)}{d_1d_2\cdots d_{\ell}}.$$\end{theorem}
\begin{proof} For simplicity we normalize in such a way that ${\rm Vol\ }S(E)=1$. We then set ${\rm Vol\ }S(\Delta)=V(R)$.  If $R$ is reducible, i.e. $R=R_1\cup R_2$ with $R_1\perp R_2$,  we have
$$V(R)=V(R_1)V(R_2).$$
Since we also have
$$\nu(R)=\nu(R_1)\nu(R_2),$$
an easy induction implies that we are reduced to show our claim under the assumption that $R$ is irreducible.

So assume $R$ irreducible and set $\alpha_0=-\theta$, with  $\theta$ the highest root. Write
$\alpha_0=\sum_{j=1}^\ell n_j\alpha_j$ with $n_j$ a negative integer  for all $j=1,\ldots ,\ell$.

As in the previous section for every $i=0,\ldots ,\ell$ set $R^{(i)}$ equal to the root system consisting of all roots in $R$ which are integral linear combinations of the roots $\alpha_0,\ldots ,\check\alpha_i ,\ldots ,\alpha_{\ell}$ so that in particular   $|R^{(i)}|\leq |R|$. Recall that the Dynkin diagram of $R^{(i)}$ is the subdiagram of $\hat D$ obtained by removing the node corresponding to $\alpha_i$. The roots $\Delta^{(i)}=\{\alpha_0,\ldots ,\check\alpha_i ,\ldots ,\alpha_{\ell}\}$ are simple roots for $R^{(i)}$. 

 We claim that   $E$ is the  union of the cones $C(\Delta^{(i)})$ whose interior are disjoint. 
To see this take $u\in E$, write $u=\sum_{h=1}^\ell b_h\alpha_h$.  If all $b_h$ are larger  or equal than zero then
$u\in C(\Delta)=C(\Delta^{(0)})$, otherwise $b_h<0$ for at least one index $1\leq h\leq \ell$. Take an index $i$ for which $b_i/n_i$ is maximum. Notice that necessarily $b_i/n_i>0$. We can clearly write
$$u=\frac{b_i}{n_i}\alpha_0+\sum_{h=1,h\neq i}^\ell (b_h-\frac{n_hb_i}{n_i})\alpha_h$$
and all coefficients are non negative.

Now observe that if, for any $i=0,\ldots ,\ell$, we write $\alpha_i$ as a linear combination of 
$\alpha_0,\ldots ,\check\alpha_i ,\ldots ,\alpha_{\ell}$ then all coefficients are negative. We then leave to the reader the easy verification that this implies that the interiors of the cones $C(\Delta^{(i)})$ are mutually disjoint. 

We deduce that
\begin{equation}\label{vol}\sum_{h=0}^\ell V(R^{(h)})=1.\end{equation}
 Now set $\Gamma=\{i|R^{(i)}=  R\}$. $\Gamma$ is not empty since $0\in \Gamma$. We can rewrite (\ref{vol}) as
$$|\Gamma|V(R)+\sum_{h\notin \Gamma}V(R^{(h)})=1.$$
Similarly by Theorem \ref{ident} we get
$$|\Gamma|\nu(R)+\sum_{h\notin \Gamma}\nu(R^{(h)})=1.$$
Since,   by the definition of $\Gamma$,   for $h\notin \Gamma$  we have $|R^{(h)}|< |R|$, by induction (the case $A_1$ in which we have 2 roots is trivial) we can assume $V(R^{(h)})=\nu(R^{(h)})$. We get
$$V(R)=\frac{1}{|\Gamma|}(1-\sum_{h\notin \Gamma}V(R^{(h)}))=\frac{1}{|\Gamma|}(1-\sum_{h\notin \Gamma}\nu(R^{(h)}))=\nu(R)$$
proving our claim
\end{proof}

\def\br#1{\langle #1\rangle}
\def\Aff#1{\widehat W}
\def\ref#1{\item{\cite{#1}}}

\newpage
\begin{appendix}

\vskip12pt

\appendix
\centerline{{\LARGE\textsc{ appendix} }}
  \medskip
  \centerline{\large {\textsc{john r. stembridge}}
  \footnote{
  Work supported by NSF grant DMS--0245385.}}

\bigskip
In this appendix, we provide an explanation for the ``curious
identity'' (Theorem~1.2)  without any case-by-case
considerations. The proof is based on two elegant formulas,
one due to L.~Solomon, the other due to R.~Steinberg.
Both of these results deserve to be better known.

If $W$ is a finite group generated by reflections in a real Euclidean
space~$E$, consider the class function on $W$ defined by
$$
\delta_W(q,t)(w):={\det(1-qw)\over\det(1-tw)}\qquad(w\in W),
$$
where the determinants are evaluated as endomorphisms of $E$,
and $q,t$ are indeterminates. This may be viewed as a bi-graded
character for $S(E)\otimes\Lambda(E)$, the tensor product of the
symmetric and exterior algebras of $E$.

In his 1963 paper on invariants of finite reflection groups~\cite{2},
Solomon explicitly determined the structure of the $W$-invariants
of $S(E)\otimes\Lambda(E)$. At the level of characters,
his structure theorem implies
$$
\br{1_W,\delta_W(q,t)}_W=\prod_{i=1}^\ell{1-qt^{d_i-1}\over 1-t^{d_i}},
\eqno(1)
$$
where $d_1,\dots,d_\ell$ are the degrees ($\ell=\dim E$),
$1_W$ denotes the trivial character of~$W$, and
$\br{f,g}_W:=|W|^{-1}\sum_{w\in W}f(w)g(w)$ is the usual
pairing of real-valued class functions $f$ and $g$.

Henceforth, assume that $W$ is a Weyl group with an irreducible
root system $R\subset E$ of rank $\ell$ and simple reflections
$S=\{s_1,\dots,s_\ell\}$.
Note that by setting $q=1$ and letting $t\to1$ in~(1),
we obtain the quantity $\nu(R)$.

We let $s_0\in W$ denote the reflection
corresponding to the highest root and set $S_0=S\cup\{s_0\}$.
One may interpret $S_0$ as the $W$-image of the simple reflections
of the associated affine Weyl group $\Aff{W}$.

Following Steinberg (see Section~3 of~\cite{3}), the action of $\Aff{W}$
on $E$ descends to a $W$-action on the $\ell$-torus $E/Q$
(where $Q$ denotes the root lattice), and the decomposition
of $E$ into simplicial alcoves by the arrangement of affine hyperplanes
associated to $R$ induces a simplicial decomposition of $E/Q$ with
a compatible $W$-action. Moreover, the $W$-stabilizers of the faces
of $E/Q$ are (up to conjugacy) generated by the various proper
subsets of $S_0$.

Given $w\in W$, Steinberg computes the Euler characteristic of the
$w$-fixed subcomplex of $E/Q$ in two different ways
(see Theorem~3.12 of~\cite{3}), thereby obtaining the identity
$$
\det(1-w) = \sum_{J\subset S_0} (-1)^{|S|-|J|} 1^W_{W_J}(w),
\eqno(2)
$$
where $W_J$ denotes the reflection subgroup generated by $J$,
and $1^W_{W_J}$ denotes the permutation character of the action
of $W$ on $W/W_J$. It is important to note that $J$ ranges over
{\it proper} subsets of $S_0$.

Steinberg actually proves a more general identity that involves
twisting by an involution; the above instance corresponds to the
trivial involution. One may also recognize~(2) as a companion to
the more familiar identity
$$
\det(w) = \sum_{J\subseteq S} (-1)^{|J|} 1^W_{W_J}(w).
$$

Now consider the evaluation of
$$
\lim_{t\to1}\ \br{\delta_W(1,0),\delta_W(1,t)}_W.
$$
First, notice that $\delta_W(1,t)\to 1_W$ as $t\to1$, so we obtain
$$
\lim_{t\to1}\ \br{\delta_W(1,0),\delta_W(1,t)}_W
  = \br{\delta_W(1,0),1_W}_W = 1
\eqno(3)
$$
by setting $(q,t)=(1,0)$ in~(1).

Second, notice that $\delta_W(1,0)(w)=\det(1-w)$, so (2) implies
$$
\br{\delta_W(1,0),\delta_W(1,t)}_W
  =\sum_{J\subset S_0}(-1)^{|S|-|J|}\br{1^W_{W_J},\delta_W(1,t)}_W
 $$$$=\sum_{J\subset S_0}(-1)^{|S|-|J|}\br{1_{W_J},\delta_{W_J}(1,t)}_{W_J},
   \eqno (4)
$$
by Frobenius reciprocity. We can evaluate each of these terms by
applying Solomon's formula to the reflection group $W_J$. But we need
to be careful, because the action of $W_J$ on $E$ will have linear
invariants if the rank of $W_J$ is less than $\ell=|S|$. In such cases,
this means that some of the degrees of $W_J$ will equal~1,
which introduces factors of $(1-q)/(1-t)$ in~(1).
Since we have set $q=1$, these factors vanish.

Thus (4) should be restricted to $\ell$-subsets of $S_0$, and we obtain
$$
\br{\delta_W(1,0),\delta_W(1,t)}_W=\sum_{j=0}^\ell\ \prod_{i=1}^\ell
  {1-t^{d_i^{(j)}-1}\over 1-t^{d_i^{(j)}}},
$$
where $d_1^{(j)}\!,\dots,d_l^{(j)}$ are the degrees of
$W_J$ for $J=S_0-\{s_j\}$.
Comparing this with~(3) in the limit $t\to1$, we obtain
the ``curious identity''
$$
\sum_{j=0}^\ell\ \prod_{i=1}^\ell 
  {d_i^{(j)}-1\over d_i^{(j)}}=1.
$$

\medskip
\end{appendix}
\vskip-\smallskipamount

\end{document}